\documentclass{article}

\usepackage{graphicx}        
\usepackage{amsfonts}
\usepackage[latin1]{inputenc}
\usepackage{latexsym,amssymb}     
\usepackage{amsmath}
\usepackage{verbatim}
\usepackage{theorem}
\usepackage{subfigure}
\usepackage{epsf}

\newtheorem{theorem}{Theorem}

\newtheorem{proposition}{Proposition}

\newtheorem{lemma}{Lemma}

\newtheorem{definition}{Definition}
\newtheorem{Res}[theorem]{Result}

\def\DEF#1{\def#1}
\renewcommand{\newcommand}{\DEF}

\def\DD{\VV^0}

\def\tribu{{\cal F}}

\newcommand{\1}{{\mathbf 1}}
\def\expect#1{\EE}

\def\invisible#1{}
\def\text#1{\quad\mbox{#1}\quad} 
\def\mtext#1{\,\mbox{#1}\,} 
\def\defegal{:=}

\def\ds{\displaystyle}
\def\eqref#1{(\ref{#1})}
\def\argmax{\arg\max} 

\newcommand{\adm}{^{\mbox \rm ad}}
\newcommand{\opt}{^{\star}}

\newcommand{\RR}{{\mathbb R}} 
 
\newcommand{\DD}{{\mathbb D}}

\newcommand{\BB}{{\mathbb B}} 
\newcommand{\EE}{{\mathbb E}} 
\newcommand{\PP}{{\mathbb P}} 
\newcommand{\UU}{{\mathbb U}} 
 
\newcommand{\VV}{{\mathbb V}} 
\newcommand{\WW}{{\mathbb W}} 
\newcommand{\XX}{{\mathbb X}}

\def\Dynamics{f}
\def\control{u}
\def\CONTROL{U}

\def\viab{{\rm Viab}}

\def\state{x}
\def\horizon{T}
\renewcommand{\omega}{w}

\def\dec{\BB} 
\def\decset#1{\dec\big(#1,\state(#1)\big)}

\def\sdo{\mathbb A} 
\def\sdoset#1{\sdo(#1)}
\def\sdoset#1{\sdo\big(#1,\omega(#1)\big)}
\def\sdosetbis#1{\sdo\big(#1,\omega\big)}

\def\unc{\mathbb S} 
\def\desir_set{\DD}

\def\feedback#1{\mathfrak{#1}}

\def\Value{\pi}
\def\viab{\VV{\rm iab}}
\def\vviab{\mbox{\tiny viab}}

\def\Unc{\Omega} 

\newcommand{\ds}{\displaystyle}

\newcommand{\VV}{\mathbb{V}}

\usepackage[dvips]{epsfig}

\begin{document}

\title{Stochastic viability and dynamic programming}
\author{ Luc \textsc{Doyen}\thanks{\textsc{CNRS, CERESP (UMR 5173
 CNRS-MNHN-P6)}, Mus\'eum National d'Histoire Naturelle, 55 rue Buffon,
 75005 Paris, France}, Michel \textsc{De Lara}\thanks{%
Universit\'e Paris--Est, \textsc{CERMICS}, 
 6--8 avenue Blaise Pascal, 77455 Marne la Vall\'ee Cedex 2, France.
}}

\maketitle

\begin{abstract}
This paper deals with the stochastic control of nonlinear systems in the
presence of state and control constraints,
for uncertain discrete-time dynamics in finite dimensional spaces. 
In the deterministic case, the viability kernel is known to play a basic
role for the analysis of such problems and the design of viable control
feedbacks.  In the present paper, we show how  a stochastic viability
kernel and viable feedbacks relying on probability (or chance) constraints can be defined and  computed by a dynamic
programming equation.   An example  illustrates most of the assertions.  
\end{abstract}

\vspace{4 mm}


    


\section{Introduction}

Risk, vulnerability, safety or precaution  constitute major issues in
the management and control of dynamical systems. Regarding these
motivations, the role played by the 
acceptability constraints or targets is central,
and it has to be articulated with 
uncertainty and, in particular, with stochasticity when a probability
distribution is given. 
The present paper addresses the issue of state and control constraints
in the stochastic context.   
For the sake of simplicity, we consider noisy control dynamics
systems. This is a natural extension of deterministic control
systems, which covers a large class of situations. 
Thus we consider the following state equation  as the uncertain dynamic model 
\begin{equation}
\state(t+1) = \Dynamics\big(t,\state(t),\control(t),
\omega(t)\big) \, , 
\quad t=t_0, \ldots, \horizon-1 \; ,
\text{ with } \state(t_{0})=\state_0
\label{eq:noisy_dyn_temps_discret}
\end{equation}
where  $\state(t)\in \XX=\RR^{n}$  represents the system \emph{state} vector at
  time $t$, $\state_0 \in \XX$ is the \emph{initial condition}
at  \emph{initial time} $t_{0}$, 
$\control(t)\in \UU=\RR^{p}$  represents
\emph{decision} or   \emph{control} vector 
while $\omega(t)\in\WW=\RR^{q}$ stands for the 
\emph{uncertain variable}\index{uncertain variable}, or
\emph{disturbance}, or \emph{noise}.

The admissibility of decisions and states is first restricted by a non empty
subset $ \dec(t,\state) $ of admissible controls in $\UU$ for all
$(t,\state)$:
\begin{equation}
\control(t)\in \decset{t}  \subset \UU  \; .
\label{eq:control_constraint_uncertain_viability}
\end{equation}
Similarly, the relevant states of the system are  limited by a non empty
subset  $\sdoset{t}$ of the  
state space $\XX$ possibly uncertain for all $ t$,
\begin{equation}
\state(t)\in \sdoset{t} \subset \XX \; ,
\label{eq:state_constraint_uncertain_viability}
\end{equation}
and a target
\begin{equation}
\state(\horizon)\in \sdoset{\horizon} \subset \XX \; .
\label{eq:target_constraint_uncertain_viability}
\end{equation}
We assume that 
\begin{equation}
 \omega(t) \in \unc(t) \subset \WW \; ,
\label{eq:unc}
\end{equation}
so that the sequences 
\begin{equation}
\omega(\cdot) \defegal \big(\omega(t_{0}), \omega(t_{0}+1), \ldots,
\omega(\horizon-1),\omega(\horizon)\big) 
\end{equation}
belonging to 
\begin{equation}
\Unc \defegal \unc(t_{0}) \times \cdots\times  \unc(\horizon) 
\subset \WW^{T+1-t_{0}}
\label{eq:Unc}
\end{equation}
capture the idea of possible \emph{scenarios} for the
problem. A scenario is an \emph{uncertainty trajectory}.  

These  control, state or target constraints may reduce the relevant paths of
the system~\eqref{eq:noisy_dyn_temps_discret}. 
Such a feasibility issue can be addressed in a robust or
stochastic framework. Here we focus on the stochastic case assuming that
the domain of scenarios $\Unc$ is equipped with some probability $\PP$.  
In this probabilistic setting, one can relax the constraint requirements
\eqref{eq:control_constraint_uncertain_viability}-%
\eqref{eq:state_constraint_uncertain_viability}-%
\eqref{eq:target_constraint_uncertain_viability}
by satisfying the state constraints along time with a given
confidence level $\beta$
\begin{equation}
	\PP \Big( \omega(\cdot) \in \Unc \mid 
	\state(t) \in \sdoset{t}	\mtext{ for } t=t_0, \ldots, T \Big) 
	\geq \beta \; 
	\label{eq:stochastic_viability_constraint}
\end{equation}
by appropriate controls satisfying 
\eqref{eq:control_constraint_uncertain_viability}. Such probabilistic constraints are often called chance constraints in the stochastic literature as in \cite{Prekopa:2003,Rockafellar:2007}.
We shall give proper mathematical content to the above formula in the
following section. Concentrating now on motivation, 
the idea of stochastic viability is basically to require the respect of
the constraints at a given confidence level $\beta$ (say $90\%$,
$99\%$). 
It implicitly assumes that some extreme events makes irrelevant the robust
approach \cite{Doyen:2000} that is closely related  to stochasticity  with
a confidence level $100\%$.  
                              
The problems of dynamic control under constraints usually refers to
viability  \cite{Aubin:1991} or invariance 
\cite{Clarkeetal:1995,Vidaletal:2000} framework.  Basically, such an approach focuses
on inter-temporal feasible paths. 
From the mathematical viewpoint, most of viability and weak invariance
results are addressed in the continuous time case. However, some
mathematical works deal with the discrete-time case. This  includes the
study of numerical schemes for the approximation of the viability
problems of the continuous dynamics as in
\cite{Aubin:1991,QSP:1995}. Important contributions for  discrete-time
case are also captured by the study of the positivity for linear systems
as in \cite{Berman-Plemmons:1994}, or by the hybrid control as in
\cite{Aubinetal:2001,Vidaletal:2000} or \cite{DeLaraetal:2006}. Other references
may be found in the control theory literature, 
such as \cite{Bitsoris:1988,Gilbert-Tan:1991} and the survey paper
\cite{Blanchini:1999}. 
A large study focusing on the discrete-time case is also provided in
\cite{DeLara-Doyen:2008}.

Viability is defined as the ability to choose,
at each time step, a control  such that the system
configuration remains admissible. 
The  \emph{viability kernel} associated
with the dynamics and the constraints  play a major role regarding such issues. It is the set of initial states $\state_0$ from which starts an acceptable solution.  For a decision maker or control designer, knowing
the viability kernel has practical interest since it describes the
states from which controls can be found that maintain the
system in a desirable configuration forever. 
However, computing this kernel is not an easy task in general. 
Of  major interest is the fact that a dynamic programming equation
underlies the computation or approximation of viability kernels as
pointed out in \cite{Aubin:1991,DeLara-Doyen:2008}. 

The present paper aims at expanding viability concepts and results in
the stochastic case for discrete-time systems.  
In particular,  we adapt the notions of viability kernel and viable
controls in the probabilistic or chance constraint framework. Mathematical materials
of stochastic viability can be found in
\cite{Aubin-Daprato:1998,Buckdahn-Quincampoix-Rainer-Rascanu:2004,Buckdahnetal:2008}
but they rather focus on the continuous time case and cope with constraints satisfied almost surely. We here provide a
dynamic programming and Bellman perspective for the probabilistic framework. 

The paper is organized as follows.
Section~\ref{sec:The_stochastic_viability_problem} is devoted to the
statement of the  probabilistic viability problem. 
Then, Section~\ref{sec:Stochastic_dynamic_programming_equation}
exhibits the dynamic programming structure underlying such stochastic
viability. 
An example is exposed in Section~\ref{sec:example} to illustrate  
some of the main findings.





\section{The stochastic viability problem} 
\label{sec:The_stochastic_viability_problem} 

Here we address the issue of state constraints in the probabilistic
sense. This is basically related to risk assessment and management.
This requires some specific tools inspired from the
viability and invariance approach known for the certain case. In particular, within the probabilistic
framework, we adapt the notions of viability kernel and viable
controls.

\subsection{Probabilistic assumptions and expected value}
\label{sec:Probabilistic_assumptions}

Probabilistic assumptions 
on  the uncertainty $\omega(\cdot) \in \Unc $ are now added,
providing a stochastic nature to the problem.
Mathematically speaking, we suppose that the domain of scenarios 
$\Unc \subset \WW^{T+1} = \RR^{q} \times \cdots \times \RR^{q} $ 
is equipped with a
$\sigma$-field\footnote{%
For instance, $\tribu$ is the trace of $\Unc$ on the usual borelian
$\sigma$-field   $\tribu = \bigotimes_{t=t_0}^{\horizon} {\cal
  B}(\RR^{q})$.} 
$\tribu$  and a \emph{probability} $\PP$: 
thus, $(\Unc,\tribu,\PP)$ constitutes a \emph{probability space}.
\index{probability!space}
The sequences
\[
\omega(\cdot)= \big(\omega(0), \omega(1), \ldots,
\omega(\horizon-1),\omega(\horizon) \big) \in \Unc  
\]
now become the \emph{primitive}\index{primitive}
\emph{random variables}.

Hereafter, 
we shall assume that 
the random process $\omega(\cdot)$ is independent and identically
distributed (i.i.d.) under probability $\PP$. 
In other words, we suppose that the probability is the product  
\( \PP = \bigotimes_{t=t_0}^{\horizon} \mu \) of a common
marginal distribution $\mu$.
The \emph{expectation operator} $\EE$ 
is defined on the set of 
measurable and integrable functions by
\[
\EE [g] 
= \EE_{\PP} \left[ g \left( \omega(\cdot) \right) \right]
= \int_{\Unc} g \big(\omega(t_0), \ldots,
\omega(\horizon) \big) 
d\mu(\omega(t_0)) \cdots d\mu(\omega(\horizon)) \; , 
\]
and we have that
\[
\EE_{\PP} \left[ g \left( \omega(t) \right) \right]
= \EE_{\mu} \left[ g \left( \omega(t) \right) \right] \; . 
\]

\subsection{Controls and feedback strategies}

It is well-known that control issues in the uncertain case are much more complicated than in the deterministic case.
In the uncertain context, we must drop the idea that the knowledge of
open-loop decisions $\control(\cdot)=
\big(\control(t_0), \ldots, \control(\horizon-1)\big)$ 
induces one single path of sequential states $\state(\cdot)=
\big(\state(t_0), \ldots, \state(\horizon)\big) $. 
Open loop controls $\control(t)$ depending only upon time $t$ are no longer
relevant, contrarily to closed loop or feedback controls 
$\feedback{\control}(t,\state(t))$ which display more adaptive properties by taking into
account the uncertain state evolution $\state(t)$. 
In the stochastic setting, all the objects considered will be 
implicitly equipped with appropriate measurability properties.
Thus we define a \emph{feedback} as an element of the set of all measurable functions
from the couples time-state towards the controls:
\begin{equation}
 \feedback{\CONTROL} \defegal 
\{ \feedback{\control} : (t,\state) \in \{t_0,\ldots,\horizon-1\} \times \XX
\mapsto \feedback{\control}(t,\state) \in \UU, \; 
\feedback{\control} \mtext{measurable}\} \; .
\label{eq:feedbackset} 
\end{equation}
The control constraints case restricts feedbacks to admissible feedbacks
accounting for control constraints
\eqref{eq:control_constraint_uncertain_viability} as follows 
\begin{equation}
\feedback{\CONTROL}\adm =\left\{\feedback{\control} \in \feedback{\CONTROL} 
\mid \feedback{\control}(t,\state)  \in \dec(t,\state) \; , \quad 
\forall (t,\state)\in \{t_0,\ldots,\horizon-1\} \times \XX \right \}  \; .
\label{eq:feedback_state_control_constraint_uncertain}
\end{equation}

Let us mention that, in the stochastic context, 
a feedback decision is also termed a  
\emph{pure Markovian strategy}. 
Markovian means that the current state contains all the
sufficient information of past system evolution to determine the
statistical distribution of future states.
Thus, only current state $\state(t)$ is needed in the feedback loop among the
whole sequence of past states $\state(t_0)$,\ldots, $\state(t)$.

At this stage, we need to introduce some notations which will appear
quite useful in the sequel:  the \emph{state map}  and 
the \emph{control map}. 
Given a feedback $\feedback{\control} \in \feedback{\CONTROL}$, 
a scenario $\omega(\cdot) \in \Unc$ 
and an initial state   $\state_0$ at time $t_0\in\{t_0,\ldots,\horizon-1 \}$,
the solution state
$\state_{\Dynamics}[t_0,\state_0,\feedback{\control},\omega(\cdot)]$ is the
state path 
$\state(\cdot)=(\state(t_0),\state(t_0+1), \ldots, \state(\horizon))$
solution of dynamics  
\begin{equation*}
\state(t+1) = \Dynamics 
\bigl(t,\state(t),\feedback{\control}(t,\state(t)),\omega(t) \bigr) 
\; , \quad t =t_0,\ldots,\horizon-1
\end{equation*}
starting from the initial condition $\state(t_0)=\state_0$ at time $t_0$ 
and associated with feedback control $\feedback{\control}$ and
scenario $\omega(\cdot)$. 
The solution control
$\control_{\Dynamics}[t_0,\state_0,\feedback{\control},\omega(\cdot)]$ 
is the associated decision path 
$\control(\cdot)=(\control(t_0),\control(t_0+1), \ldots, \control(\horizon-1))$ 
where $\control(t)=\feedback{\control}(t,\state(t))$.

\subsection{The stochastic viability kernel and viable feedbacks}

The viability kernel plays a major role in the viability analysis. In
the determinsitic case, it is the set of initial states $\state_0$ such that
the state constraints hold true for at least one control stategy.  
In the probabilistic setting, one relaxes the constraints requirement by
satisfying the state constraints along time with a given 
confidence level as in \eqref{eq:stochastic_viability_constraint}.
We give proper mathematical content to this latter
formula~\eqref{eq:stochastic_viability_constraint} inspired by chance constraints \cite{Prekopa:2003} in the following
Definition. 

\begin{definition} 
The \emph{stochastic viability kernel} at time $t_0$ and at confidence level
$\beta\in ]0,1]$ is 
\index{stochastic!viability kernel} 
\index{viability kernel!stochastic}
\begin{equation}
 \viab_{\beta}(t_0) \defegal \left\{ \state_0 \in \XX   \;\left|
\begin{array}{l} \mtext{ there exists } \feedback{\control} \in \feedback{\CONTROL}\adm 
\mtext{ such that}\\
\PP \Big(\omega(\cdot) \in \Unc \mid 
\state(t) \in \sdoset{t} \mtext{ for } t=t_0, \ldots, T \Big) \geq \beta 
\end{array} \right. \right\}
\label{eq:stochastic_viability_kernel}
\end{equation}
\noindent where $\state(t)$ 
is a shorthand for the solution map
$\state(t)=\state_{\Dynamics}[t_0,\state_0,\feedback{\control},\omega(\cdot)](t)$.
\end{definition}

Stochastic viable feedbacks are measurable feedback controls that allow the
stochastic viability property to hold true. 

\begin{definition}
\emph{Stochastic viable feedbacks} are those  for which
the above relations occur:
\begin{equation}
\feedback{\CONTROL}_{\beta}^{\vviab}(t_0,\state_0) \defegal \left\{\feedback{\control}\in \feedback{\CONTROL}\adm  \left|\; 
\begin{array}{l} \PP \Big(\omega(\cdot) \in \Unc \mid \state(t) \in \sdoset{t} 
\mtext{ for } t=t_0, \ldots, T \Big) \geq \beta 
\end{array} \right. \right\} . 
\label{eq:stochastic_viable_feedbacks}
\end{equation}
\end{definition}

We have the following strong link between stochastic viable feedbacks
and the viability kernel:  
\begin{equation*}
\state_0 \in \viab_{\beta}(t_0)  \iff 
\feedback{\CONTROL}_{\beta}^{\vviab}(t_0,\state_0) \not = \emptyset \; .
\end{equation*}

Of particular interest is the case where the confident rate is $\beta=1$
which is very close to robust viability and control. 
Indeed, when the scenario domain $\Unc$ is countable and that
every scenario $\omega(\cdot)$ has strictly positive probability under
$\PP$, \( \viab_{1}(t_0) \) is the \emph{robust viability kernel}
(the set of initial states $\state_0$ such that
the state constraints hold true for at least one control stategy,
\emph{whatever the scenario}).
When the uncertainty domain $\unc(t)$ in~\eqref{eq:unc} is reduced to a
single element, so is also the scenario domain $\Unc$ in~\eqref{eq:Unc}:
this is the deterministic case for which \( \viab_{1}(t_0) \) coincides
with the classical viability kernel
\cite{Aubin:1991,DeLara-Doyen:2008}. 

\section{Stochastic dynamic programming equation}
\label{sec:Stochastic_dynamic_programming_equation}

We shall now exhibit a a characterization  of stochastic viability in
terms of  dynamic programming.
It relies on the  the \emph{maximal viability probability} defined recursively
as follows. 

\begin{definition} 
Assume that the random process $\omega(\cdot)$ is i.i.d. under
probability $\PP$, with marginal distribution \( \mu \). 
The \emph{stochastic viability value function} $V(t,\state)$, 
associated with dynamics~\eqref{eq:noisy_dyn_temps_discret}, 
control constraints~\eqref{eq:control_constraint_uncertain_viability}, state constraints~\eqref{eq:state_constraint_uncertain_viability}
and target constraints~\eqref{eq:target_constraint_uncertain_viability}
is defined by the following backward induction:
\begin{equation}
\left\{ \begin{array}{rcl}
V(\horizon,\state) & \defegal & \EE_{\mu}
\biggl[\1_{ \sdosetbis{\horizon} }(\state)\biggr] ,\\[4mm]
V(t,\state) & \defegal & \displaystyle 
 \max_{\control \in \dec(t,\state)} 
\EE_{\mu} \biggl[ \1_{ \sdosetbis{t} }(\state) \; 
V\big(t+1,\Dynamics(t,\state,\control,\omega) \big) \biggr] \; .
\end{array} \right.
\label{eq:Bellman_discret_stochastic_viability}
\end{equation}
\end{definition}
Here, $\1_A$ stands for the indicator function  of a set $A$.
It is defined by $\1_A(\state)=1$ if $\state \in A$, 
and $\1_A(\state)=0$ if $\state \not \in A$.


The backward dynamic 
programming equation~(\ref{eq:Bellman_discret_stochastic_viability})  
allows us to define the value function $V(t,\state)$. By writting a $\max$ instead of a $\sup$, we implicitly assume the existence of an optimal solution for each time $t$ and state $x$.  
It turns out that the stochastic viability function $V(t_{0},\state)$
at time $t_{0}$
is related to the stochastic viability kernels 
\( \{ \viab_{\beta}(t_0) , \beta \in [0,1] \} \), and that
dynamic programming induction reveals relevant stochastic feedback
controls. 
To achieve this,  we first claim that the value function $V$ is the
solution of a (stochastic) optimal control problem involving the
viability criterion $\Value$ defined as follows: 
\begin{equation}
	\Value\bigl(t_{0},\state(\cdot),\control(\cdot),\omega(\cdot)\bigr) = 
		\prod_{t=t_{0}}^{\horizon} \1_{ \sdoset{t} }(\state(t)) \; .
	\label{eq:value}
\end{equation}

\begin{proposition}
Assume that the random process $\omega(\cdot)$ is i.i.d. under
probability $\PP$, with marginal distribution \( \mu \). 
For any initial conditions $(t_0,\state_0)$, we have
\[
V(t_0,\state_0)=
\max_{\feedback{\control}\in \feedback{\CONTROL}\adm} \EE_{\PP}
\left[\Value\bigl(t_{0},\state(\cdot),\control(\cdot),\omega(\cdot)\bigr) \right]\; ,
\]
where the stochastic viability value function \( V(t_0,\state_0) \) is
given by the backward
induction~\eqref{eq:Bellman_discret_stochastic_viability},
where the criterion $\Value$ is defined in \eqref{eq:value}, and   
where
$\state(\cdot)=\state_{\Dynamics}
[t_0,\state_0,\feedback{\control},\omega(\cdot)](\cdot)$
and 
$\control(\cdot)=
\control_{\Dynamics}[t_0,\state_0,\feedback{\control},\omega(\cdot)]$ 
are shorthand expressions for the solution maps.
\label{prop:opti_viability_proba}
\end{proposition}

The proof of this previous Proposition is exposed in Appendix \ref{sec:proofs}. 
We also derive the following assertion regarding the 
stochastic viability kernel.

\begin{proposition} 
Assume that the random process $\omega(\cdot)$ is i.i.d. under
probability $\PP$, with marginal distribution \( \mu \). 
The stochastic viability kernel at confidence level $\beta$ 
is the section of level $\beta$ of the stochastic value function:  
\begin{equation}
V(t_0,\state_0) \geq \beta \iff \state_0 \in \viab_{\beta}(t_0) \; .
\end{equation}
\label{pr:Bellman_discret_stochastic_viability_noyau}
\end{proposition}

The proof of this previous Proposition is also exposed in Appendix \ref{sec:proofs}. 
As regard the viable feedbacks, we obtain the following assertion.

\begin{proposition} 
Assume that the random process $\omega(\cdot)$ is i.i.d. under
probability $\PP$, with marginal distribution \( \mu \).   
For any time $t=t_{0},\ldots,\horizon-1$ and state $\state$, let us assume that
\begin{equation}
\dec^{\vviab} (t,\state) \defegal \argmax_{\control \in \dec(t,\state)}
\EE_{\mu} \Biggl[ \1_{ \sdosetbis{t} }(\state)  \; 
V \biggl( t+1, \Dynamics \bigl( t,\state,\control,\omega(t) \bigr)
\biggr) \Biggr]  
\end{equation}
is not empty.
Then,  for any $\state_0\in\viab_{\beta}(t_0)$, any measurable
selection\footnote{Any  $\feedback{\control}\opt \in\feedback{\CONTROL}$
such that $ \feedback{\control}\opt(t,\state) \in \dec^{\vviab} (t,\state)$ for any $t$ and $x$.}
$\feedback{\control}\opt \in \dec^{\vviab}$ belongs to the set of
stochastic viable feedbacks  
$\feedback{\CONTROL}_{\beta}^{\vviab}(t_0,\state_0)$. 
\label{pr:Bellman_discret_stochastic_viability_control}
\end{proposition}


\section{A simple academic example}
\label{sec:example}

  To illustrate the general statements, we consider a simple academic
  model and perform a probabilistic viability analysis. 

\subsection{Example statement}
 
The evolution of a scalar  $\state(t)$ is governed by the discrete-time
dynamics 
\[\state(t+1)=\state(t) + \control(t) + \omega(t) \; ,
\]
where control is constrained by 
\[
\control(t)\in\{-1,1\}=\dec(t,\state)=\dec 
\]
and uncertainty scenarii are induced by 
\[
\omega(t)\in \{-1,0,1\}=\unc(t)=\unc \; .
\]
We assume that $\omega(\cdot)$ is an i.i.d. sequence, 
with probability 
\[
\mu(\omega(t)=1)= \mu(\omega(t)=-1)=p; \;\mu(\omega(t)=0)=1-2p \; .
\]
The state constraint is
\[
	\state(t)\in \{-1,0,1\} =\sdoset{t}=  \sdo\; . 
\]
The decision maker  intends to exhibit controls such that this
constraint is satisfied with a high enough probability 
\[
\PP\biggl(\state(t)\in \{-1,0,1\} \; , 
\quad t=t_0,\ldots,\horizon \biggr) \geq \beta \; .
\]

The intuition to satisfy the above probability constraint is as
follows. When $\state(t)$ belongs to the border \( \{-1,1\} \) 
of the domain \( \sdo = \{-1,0,1\} \), there is an obvious decision to
make: if $\state(t)=-1$, take \( \control(t) =1 \) so that 
\( \state(t) + \control(t) = 0 \) and thus
\( \state(t+1)= \omega(t)  \in \{-1,0,1\} \) (the same with 
$\state(t)=1$ and \( \control(t) =-1 \) ).
But when $\state(t)=0$, then \( \state(t+1)= \control(t) + \omega(t) \) 
and, whatever \( \control(t) \in\{-1,1\} \), there is a chance that 
$ \omega(t) $ takes the same value, sending $\state(t)$ outside 
\( \sdo = \{-1,0,1\} \).

\subsection{Results}

\begin{figure}[ht]	
\centering
\epsfxsize=12 cm\epsfbox{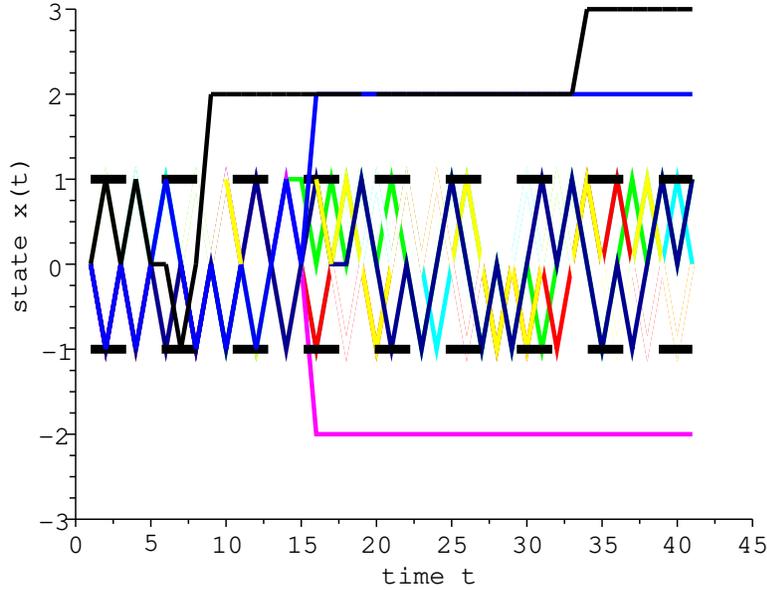}
\caption{\label{fig:viab} 9 simulations of state trajectories $x(t)$ over time horizon $[0,40]$ for dynamics $x(t+1)=x(t)+u(t)+w(t)$ starting from $x_0=0$ with stochastic viable feedback controls $\feedback{\control}\opt(t,x)\in \dec^{\vviab} (t,\state)$ as defined in (\ref{eq:decviabexample}). Probability  of facing high disturbances $w\in\{-1,+1\}$ is low with $p=1\%$. Viability probability value function $V(0,0)\approx 67\%$ and 3 trajectories over 9 violate the constraint.} \end{figure}	

By dynamic programming
equation~(\ref{eq:Bellman_discret_stochastic_viability}), we compute the
maximal viability probability 
$V(t,\state)$ and associated viable feedback
controls $\dec^{\vviab} (t,\state)$. 

\begin{Res} 
Introduce  matrix $M$, vectors $\vec{\1}$ and
 $\vec{\1_{i}}(\state)$ by
 \[
M= \left( \begin{array}{ccc}
p & 1-2p & p \\
p &1-2p &0 \\
p &1-2p &p
\end{array} \right), \;
\vec{\1}=\left( \begin{array}{c}
1  \\
1  \\
1 
\end{array} \right), \;
(\vec{\1}_i(\state))_j=\1_{\{i\}}(\state)=\1_{\{\state=i\}}   \; . 
\]
The stochastic viability value function is given by
\[
V(t,\state)=\sum_{i=-1,0,1}\left<\vec{\1_{i}}(\state), 
M^{\horizon-t}\vec{\1}\right> \; ,
\]
or, in other words, \( V(t,\state)= 0 \) for all 
\(  \state \not\in \{-1,0,1\} \) and 
\[
V(t,\state)=(M^{\horizon-t}\vec{\1})_{\state+2} \; , \quad
\forall \state \in \{-1,0,1\} \; .
\]
The associated viable feedback controls are given by

\begin{equation}
	\dec^{\vviab} (t,\state)= \left\{ \begin{array}{lcl}
	1 & \mtext{ if } & \state=-1 \\
	\{-1,1\} & \mtext{ if } & \state=0 \\
	-1 & \mtext{ if } & \state=1 \; ,\\
	\end{array} \right.
	\label{eq:decviabexample}
\end{equation}

\end{Res}

Consequently the viability kernel reads:
\begin{Res} 
	\[\viab_\beta(t)=\left\{\begin{array}{ccc}
	\sdo & \mbox{ if }& \beta \leq  (M^{\horizon-t}\vec{\1})_{2}\\
	\{-1,1\} & \mbox{ if }&    (M^{\horizon-t}\vec{\1})_{2} <  \beta \leq  (M^{\horizon-t}\vec{\1})_{1}\\
	\emptyset & \mbox{ if }&     (M^{\horizon-t}\vec{\1})_{1} < \beta. 
	\end{array}\right.\]
\end{Res}

The difficulty of the control is captured by  the second row of the matrix $M$ where the sum is not equal to 1 which suggests that the state $x=0$ can escape from $\sdo$.    The results  are illustrated by Figure \ref{fig:viab} where 9 simulations of state trajectories $x(t)$  starting from $x_0=0$ are displayed over time horizon $[0,40]$ with stochastic viable feedback controls $\feedback{\control}\opt(t,x)\in \dec^{\vviab} (t,\state)$ as defined in (\ref{eq:decviabexample}). Probability  of facing high disturbances $w\in\{-1,+1\}$ is low with $p=1\%$. However viability probability value function turns out to be  $V(0,0)\approx 67\%$ which points out  a significant risk of leaving viability set $\sdo=\{-1,0,1\}$ due the accumulation of risks over 40 periods; Therefore it is  intuitive that 3 paths over 9 leave the state constraint set $\sdo=\{-1,0,1\}$ along time.


\vspace{2 mm}

\begin{proof}
We shall check that \( V(t,\state)=\sum_{i=-1,0,1}\left<\vec{\1_{i}}(\state), 
M^{\horizon-t}\vec{\1}\right> \) is solution to the dynamic programming
equation~(\ref{eq:Bellman_discret_stochastic_viability}).

This is true for final time $t=\horizon$ because
\[
\sum_{i=-1,0,1}\left<\vec{\1_{i}}(\state), M^{\horizon-\horizon}\vec{\1}\right>
=\sum_{i=-1,0,1}\left<\vec{\1_{i}}(\state),\vec{\1}\right>
=\sum_{i=-1,0,1}\1_{i}(\state)=\1_{\{ -1,0,1 \} }(\state) 
=\1_{\sdo}(\state) \; . 
\]
Proceeding by backward induction, let us suppose that
\[
V(t+1,\state)=
\sum_{i=-1,0,1}\left<\vec{\1_{i}}(\state),
  M^{\horizon-(t+1)}\vec{\1}\right>  \; . 
\]
The dynamic programming
equation~(\ref{eq:Bellman_discret_stochastic_viability}) gives
\[
V(t,\state) = \displaystyle 
\1_{ \{-1,0,1\} }(\state) \; 
\max_{\control \in\{-1,1\} } 
\EE_{\mu} \biggl[ 
V\big(t+1,\state + \control + \omega \big) \biggr] \; .
\]
Whenever $\state\notin\sdo=\{-1,0,1\} $, we clearly have that 
$V(t,\state)=0$.  
 Whenever $\state=-1$, we deduce 
that   
\[
	\begin{array}{rcl}
V(t,-1)		&=&\ds \max\bigl\{
pV(t+1,-3)+ (1-2p) V(t+1,-2) + pV(t+1,-1), \\
&& pV(t+1,-1)+ (1-2p) V(t+1,0) + pV(t+1,1) \bigr\} \\
		&=&\ds \max\bigl\{ pV(t+1,-1)+  (1-2p)V(t+1,0) + pV(t+1,1), 
pV(t+1,-1) \bigr\} \\
		&=&\ds pV(t+1,-1)+  (1-2p)V(t+1,0) + pV(t+1,1)
\end{array}
\]
and the viable control is provided by
	\( \feedback{\control}\opt(t,-1)=1.\)
By induction, we deduce that
	\[
	\begin{array}{rcl}
V(t,-1) &=&\ds pV(t+1,-1)+  (1-2p)V(t+1,0) + pV(t+1,1)\\
 &=&\ds \sum_{i=-1,0,1}M_{1,i+2}(M^{\horizon-(t+1)}\vec{\1})_{i+2}\\
  &=&\ds (MM^{\horizon-(t+1)}\vec{\1})_{1}\\
    &=&\ds (M^{\horizon-t}\vec{\1})_{1}\\
 &=&\ds \left<\vec{\1_{-1}}(-1), M^{\horizon-t}\vec{\1}\right>\\
 &=&\ds \sum_{i=-1,0,1}\left<\vec{\1_{i}}(-1),
   M^{\horizon-t}\vec{\1}\right> \; . 
\end{array}
\]

In the same way, we check the expression for the 
stochastic viability value function $ V(t,1) $ when $\state=1$,
and obtain the viable control \( \feedback{\control}\opt(t,1)=-1 \).
The case $\state=0$ is treated in the same vein, with the difference
that viable control is not unique since 
	\( \feedback{\control}\opt(t,0)\in\{-1,+1\} \)
and 
	\[
V(t,0) =\ds pV(t+1,-1)+  (1-2p)V(t+1,0).
\]

\end{proof}

\newpage
\appendix

\section{Proofs}\label{sec:proofs}

\subsection{Proof of Proposition \ref{prop:opti_viability_proba}}


We use the following notations for any strategy 
$\feedback{\control}\in \feedback{\CONTROL}$:

\begin{itemize}
	\item  $\Value^{\feedback{\control}}$  is the evaluation of the
      criterion $\Value$ defined in \eqref{eq:value} 
\begin{equation}
\Value^{\feedback{\control}}\bigl(t_0,\state_0,\omega(\cdot)\bigr) \defegal
\Value\bigl(t_{0},\state_{\Dynamics}[t_0,\state_0,\feedback{\control},\omega(\cdot)](\cdot),
\control_{\Dynamics}[t_0,\state_0,\feedback{\control},\omega(\cdot)](\cdot),\omega(\cdot)\bigr)
\label{eq:value_feedback_uncertain_optimality_proof}
\end{equation}
\noindent where 
$\omega(\cdot) \in \Unc$ and $\state_{\Dynamics}$, 
$\control_{\Dynamics}$ are the solution maps;

\item the expected value 
\begin{equation}
 \Value^{\feedback{\control}}_{\EE}(t_0,\state_0) \defegal 
\EE_{\PP} 
\left[
  \Value^{\feedback{\control}}\bigl(t_0,\state_0,\omega(\cdot)\bigr)
\right] \, . 
\label{eq:GG_performance}
\end{equation}
\end{itemize}

We consider the maximization problem:
\begin{equation}
 \Value_{\EE}\opt(t_{0},\state_0) \defegal \displaystyle
\max_{\feedback{\control}\in \feedback{\CONTROL}\adm} 
\Value^{\feedback{\control}}_{\EE}(t_0,\state_0) \, .
\label{eq:value_payoff_discret_Whittle}
\end{equation}
We  aim at proving that 
\[
V(t,\state)=\Value^{\feedback{\control}\opt}_{\EE}(t,\state)= 
\max_{\feedback{\control}\in \feedback{\CONTROL}\adm} \Value^{\feedback{\control}}_{\EE}(t,\state) \; .
\]
Let $\feedback{\control}\opt \in \feedback{\CONTROL}\adm$ denote one of
the measurable viable feedback strategies given by the dynamic
programming equation~\eqref{eq:Bellman_discret_stochastic_viability}. 
We perform a backward induction to
prove~\eqref{eq:value_payoff_discret_Whittle}. 

First, the equality at $t=\horizon$ holds true since
\[
\begin{array}{rcll}
\Value^{\feedback{\control}\opt}_{\EE}(\horizon,\state) &=& \ds \EE_{\PP} 
\left[ \Value^{\feedback{\control}\opt}(\horizon,\state,\omega(\cdot)) \right] 
& \text{by definition ~\eqref{eq:GG_performance}} \\
&=& \ds \EE_{\mu} \left[ \1_{\sdosetbis{\horizon}}(\state) \right] 
& \text{by definition ~\eqref{eq:value}} \\
&=& V(\horizon,\state) & \text{by definition ~\eqref{eq:Bellman_discret_stochastic_viability}.} 
\end{array}
\]
Now, suppose that 
\begin{equation}
\Value^{\feedback{\control}\opt}_{\EE}(t+1,x)=
\max_{\feedback{\control}\in \feedback{\CONTROL}\adm} \Value^{\feedback{\control}}_{\EE}(t+1,x)=V(t+1,x)
\; .
\label{eq:induction_hypothesis}
\end{equation}
The very definition (\ref{eq:Bellman_discret_stochastic_viability}) 
of the value function $V$ by dynamic programming combined with~\eqref{eq:preparation_Bellman} in
Lemma~\ref{lem-preparation_Bellman} (proved below) imply that 
\[
\begin{array}{rcll}
\Value^{\feedback{\control}\opt}_{\EE}(t,\state) &=&\EE_{\mu} \left[ \1_{\sdoset{t}}(\state)\; \Value^{\feedback{\control}\opt}_{\EE}\left(t+1,\Dynamics\big(t,\state,\feedback{\control}\opt(t,\state),\omega(t)) \right)  \right] 
& \text{by~\eqref{eq:preparation_Bellman}} \\
&= &\EE_{\mu} \left[ \1_{\sdosetbis{t}}(\state) \;V\left(t+1,\Dynamics\big(t,\state,\feedback{\control}\opt(t,\state),\omega) \right)  \right] 
& \text{by~\eqref{eq:induction_hypothesis}} \\
&= &\max_{\control \in \dec(t,\state)} \EE_{\mu}  \left[ \1_{\sdosetbis{t}}(\state) \; V\left(t+1,\Dynamics\big(t,\state,\control,\omega)\right) \right] 
& \text{by~\eqref{eq:Bellman_discret_stochastic_viability}} \\
& = & V(t,\state) 
& \text{by~\eqref{eq:Bellman_discret_stochastic_viability}.} 
\end{array}
\]

Similarly, for any $\feedback{\control}\in\feedback{\CONTROL}\adm$, we obtain 
\[
\begin{array}{rcll}
\Value^{\feedback{\control}}_{\EE}(t,\state) &=&\EE_{\mu} \left[ \1_{\sdoset{t}}(\state) \;\Value^{\feedback{\control}}_{\EE}\left(t+1,\Dynamics\big(t,\state,\feedback{\control}(t,\state),\omega(t)) \right) \right] 
& \text{by~\eqref{eq:preparation_Bellman}} \\
&\leq  &\EE_{\mu} \left[ \1_{\sdosetbis{t}}(\state) \; V\left(t+1,\Dynamics\big(t,\state,\feedback{\control}(t,\state),\omega) \right)   \right] 
& \text{by~\eqref{eq:induction_hypothesis}} \\
& \leq &\max_{\control \in \dec(t,\state)}\EE_{\mu} \left[ \1_{\sdosetbis{t}}(\state) \; V\left(t+1,\Dynamics\big(t,\state,\control,\omega)\right) \right] 
& \mtext{ since } \feedback{\control}(t,\state) \in \dec(t,\state) \\
& = & V(t,\state) & \text{by~\eqref{eq:Bellman_discret_stochastic_viability}.} 
\end{array}
\]
Consequently, the desired statement is obtained  since 
\[
\max_{\feedback{\control}\in \feedback{\CONTROL}\adm} \Value^{\feedback{\control}}_{\EE}(t,\state)\leq V(t,\state)  =\Value^{\feedback{\control}\opt}_{\EE}(t,\state)
\] yields the equality
\[
V(t,\state)=\Value^{\feedback{\control}\opt}_{\EE}(t,\state)= 
\max_{\feedback{\control}\in \feedback{\CONTROL}\adm} \Value^{\feedback{\control}}_{\EE}(t,\state) \; .
\]

\begin{lemma}
We have, for  $t=t_{0},\ldots,\horizon-1$ and $\feedback{\control} \in \feedback{\CONTROL}$,
\begin{equation}
\left\{ \begin{array}{rl}
\Value^{\feedback{\control}}_{\EE}(T,\state) =& 
\EE_{\mu} \left[ \1_{\sdoset{T}}(\state) \right] \\[4mm]
\Value^{\feedback{\control}}_{\EE}(t,\state) =&\EE_{\mu}  \left[ 
\1_{\sdosetbis{t}}(\state) \; \Value^{\feedback{\control}}_{\EE}\left(t+1,\Dynamics\big(t,\state,\feedback{\control}(t,\state),\omega) 
\right)   \right] \; .
\end{array} \right.
\label{eq:preparation_Bellman}
\end{equation}
\label{lem-preparation_Bellman}
\end{lemma}

\begin{proof}
By~\eqref{eq:value} and
\eqref{eq:value_feedback_uncertain_optimality_proof}, we have 
\begin{equation}
\left\{ \begin{array}{rcl}
\Value^{\feedback{\control}}(T,\state,\omega(\cdot)) &=& \1_{\sdoset{T}}(\state) \\[4mm] 
\Value^{\feedback{\control}}(t,\state,\omega(\cdot)) &=& \1_{\sdoset{t}}(\state)
 \Value^{\feedback{\control}}(t+1,\Dynamics\big(t,\state,\feedback{\control}(t,\state),\omega(t)),\omega(\cdot)) ) \; .
\end{array} \right.
\label{eq:P(t,\state,gamma,\omega)}
\end{equation}
Notice that 
$\Value^{\feedback{\control}}(t,\state,\omega(\cdot))$ depends only upon 
the end $ \big(\omega(t), \ldots, \omega(\horizon-1) \big)$,
and not upon 
the beginning $\big(\omega(t_0), \ldots, \omega(t-1) \big)$.
We shall write this property abusively by
\begin{equation}
\Value^{\feedback{\control}}(t,\state,\omega(\cdot)) =
\Value^{\feedback{\control}}
\left(t,\state,\left(\omega(t), \ldots, \omega(\horizon-1)\right) \right) \; .
\label{eq:causality_P(t,\state,gamma,\omega)}
\end{equation}
We have
\[
\begin{array}{rcll}
\Value^{\feedback{\control}}_{\EE}(t,\state) &=& \ds \EE_{\PP}
\left[ \Value^{\feedback{\control}}(t,\state,\omega(\cdot)) \right]
 &\text{by~\eqref{eq:GG_performance}} \\
 &=& \ds \EE_{\PP} \left[ \1_{\sdoset{t}}(\state)\; \Value^{\feedback{\control}}\left(t+1,\Dynamics\big(t,\state,\feedback{\control}(t,\state),\omega(t)),\omega(\cdot) \right)   \right] & \text{by~\eqref{eq:P(t,\state,gamma,\omega)}} \\
 &=& \displaystyle \EE_{\mu}  
\biggl[ \EE_{\mu^{\horizon-t-1}} \bigl[  \bigr. \bigr. &\\
& & \hspace{-1 cm}\biggl. \bigl. \1_{\sdoset{t}}(\state)
\Value^{\feedback{\control}}\left(t+1,\Dynamics\big(t,\state,\feedback{\control}(t,\state),\omega(t)),\omega(t+1),\ldots,\omega(\horizon-1)
\right)  \bigr]  \biggr]& \\
&&\mtext{by Fubini theorem} & \\ 
&=& \ds \EE_{\mu} \biggl[\1_{\sdoset{t}}(\state)  \biggr. & \\
&&   \hspace{-1 cm} \biggl.\; \EE_{\mu^{\horizon-t-1}} 
\left[\Value^{\feedback{\control}}\left(t+1,\Dynamics\big(t,\state,\feedback{\control}(t,\state),\omega(t+1),\ldots,\omega(\horizon-1)) \right)\right] \biggr] & \\
&=& \ds \EE_{\mu}  \biggl[\1_{\sdoset{t}}(\state)  \EE_{\omega(\cdot)\in\Unc} \bigl[\Value^{\feedback{\control}}\left(t+1,F\left(t,\state,\feedback{\control}(t,\state),\omega(\cdot)\right)\right) \bigr]    \biggr] 
& \text{by~\eqref{eq:causality_P(t,\state,gamma,w)}}  \\ 
&=& \ds \EE_{\mu} \bigl[\1_{\sdosetbis{t}}(\state)
\;\Value^{\feedback{\control}}_{\EE}(t+1,\Dynamics\big(t,\state,\feedback{\control}(t,\state),\omega \big)) \bigr] 
& \text{by~\eqref{eq:GG_performance}.} 
\end{array}
\]

\end{proof}


\subsection*{Proof of Proposition~\ref{pr:Bellman_discret_stochastic_viability_noyau}} 

It is enough to remark  that
\begin{equation}
\viab_{\beta}(t) = \left\{\state_0\in \XX \;\left|
\max_{\control(\cdot)} \EE_{\PP} \left[\Value\bigl(t_{0},\state(\cdot),\control(\cdot),\omega(\cdot)\bigr)\right]
 \geq  \beta
\right.\right\} \; .
\end{equation}


\subsection*{Proof of Proposition~\ref{pr:Bellman_discret_stochastic_viability_control}} 

Simply follow step by step the proof of Proposition
\ref{prop:opti_viability_proba}.


\newcommand{\noopsort}[1]{} \ifx\undefined\allcaps\def\allcaps#1{#1}\fi

\end{document}